\input amssym.def
\input amssym
\documentstyle[12pt]{article}
\newcommand{\less}{\setminus}

\newcommand{\cl}[1]{\overline{#1}}
\newcommand{\Cal}[1]{{\cal#1}}

\newcommand\Fn{\mbox{Fn}}
\newcommand\forces{\Vdash}

\newcommand\from{\!:\!}
\newcommand\re{\!\restriction\!}
\newcommand\bbb[1]{{\Bbb#1}}
\newcommand\suppt{\mbox{suppt}}
\newcommand\dom{\mbox{dom}}
\newcommand\cov{\mbox{cov}}
\newtheorem{Theorem}{Theorem}[section]
\newtheorem{Lemma}[Theorem]{Lemma}

\newtheorem{Def}[Theorem]{Definition}
\newenvironment{Definition}{\begin{Def}\em}{\end{Def}}
\newenvironment{Proof}{\noindent{\em Proof.}}{~$\Box$\bigskip}

\title{Forcing tightness in products of fans} \author{ J\"org 
Brendle\footnotemark\ 
and Tim LaBerge}
\begin{document}
\addtolength{\baselineskip}{.5\baselineskip}

\maketitle
\begin{abstract} 
We define several notions of forcing that allow us to manipulate the 
tightness of products of fans. Some consequences include: 
$t(F_\theta\times F_\omega) = \theta$ does not imply  the 
existence of a 
$(\theta,\omega)$-gap,  new  examples of   first countable 
$<\!\theta$-cwH spaces that are not $\le\!\theta$-cwH for   
singular 
cardinals $\theta$, and for cardinals $\lambda\le\theta$ with 
$cf(\theta)\ge \omega_1$ and either $\lambda$ regular or 
$\lambda^\omega\le\theta$, a first 
countable 
$<\!\theta$-cwH not $\le\!\theta$-cwH space that can be made cwH 
by removing a closed discrete set of cardinality $\lambda$. 
We also prove two theorems that characterize tightness of products 
of fans in terms of families of integer-valued functions.

\end{abstract}

\renewcommand{\thefootnote}{\fnsymbol{footnote}}
\footnotetext[1]{Supported by DFG--grant Nr. Br 1420/1--1.} 

\section{Introduction}  The $\theta$-fan $F_\theta$ is the quotient 
space obtained by 
identifying the non-isolated points of the product 
$\theta\times(\omega+1)$ to a single point $\infty$. (Here $\theta$ 
has the discrete topology and $\omega+1$ has the order topology.) 
Thus, a neighborhood of $\infty$ is a set of the form
$$V_g = \{\infty\}\cup\{\langle\alpha,m\rangle :m> 
g(\alpha)\}\qquad g\in\omega^\theta.
$$
When $\lambda\le \theta$, we use sets 
$$V_g \times U_f= 
\{(\infty,\infty)\}\cup\{(\langle\beta,m\rangle,\langle\alpha,n
\rangle): m> g(\beta) \land n> f(\alpha)\}\qquad 
g\in\omega^\theta\mbox{, }f\in\omega^\lambda
$$
as a base at $(\infty,\infty)$ in the product $F_\theta\times 
F_\lambda$.

The tightness $t(p,X)$ of a point $p$ in a topological space $X$ is the 
supremum of the cardinalities of all $A\subseteq X$ such that $p\in 
\cl{A}$, but whenever $B\subseteq A$ and $|B|<|A|$, then 
$p\notin\cl{B}$. The tightness of $X$ is then $t(X) = \sup\{t(p,X):p\in 
X\}$. The tightness of $X$ is the least upper bound of the 
cardinalities of the subsets of $X$ needed to define the closure 
operator.

Clearly, the tightness of $F_\theta$ is $\omega$, and it is not hard 
to see that $t(F_\theta\times F_\lambda) = 
t((\infty,\infty),F_\theta\times F_\lambda)$ (details are in 
\cite{LL94}). This motivates the following 
definition.

\begin{Definition} Let $\lambda\le\theta$ be infinite cardinals and 
suppose that $A\subseteq(\theta\times 
\omega)\times(\lambda\times \omega)$.
\begin{enumerate}
\item[(1)] If $\lambda<\theta$, we say that $A$ is 
$(\theta,\lambda)$-good if
\begin{enumerate}
\item[(a)] $(\infty,\infty)\in\cl{A}$;
\item[(b)] $\forall B\in A^{<\theta}\;((\infty,\infty)\notin \cl{B})$; and
\item[(c)] $\forall 
E\in[\lambda]^{<\lambda}\;((\infty,\infty)\notin\cl{A\cap 
((\theta\times\omega)\times (E\times 
\omega))}\,)$.
\end{enumerate}
\item[(2)] If $\lambda=\theta$, we say that $A$ is 
$(\theta,\theta)$-good if (a), (b),  and the following are true:
\begin{enumerate}
 \item[(c)] $\forall 
E\in[\theta]^{<\theta}\;((\infty,\infty)\notin\cl{A\cap 
((\theta\times\omega)\times (E\times 
\omega))}\,)$ and
\item[(d)] $\forall F\in 
[\theta]^{<\theta}\;((\infty,\infty)\notin\cl{A\cap ((F\times 
\omega)\times 
(\theta\times\omega))}\,)$.
\end{enumerate}
\end{enumerate}
\end{Definition}

The existence of a $(\theta,\lambda)$-good set $A$ implies that 
$t(F_\theta\times F_\lambda) = \theta$; moreover, if either 
$\theta'<\theta$ and $\lambda'\le\lambda$ or $\theta'\le\theta$ and 
$\lambda'<\lambda$, then $A$ cannot be construed as a subset of 
$F_{\theta'}\times F_{\lambda'}$.

As shown in \cite{LL94}, the 
existence of a $(\theta,\lambda)$-good set  is equivalent to the 
existence of a first countable 
$<\!\theta$-cwH space $X$ with a closed discrete set $D$ of 
cardinality $\theta$ such that $D$ is not separated and 
$$\lambda = 
\min\{|E|:E\subseteq D\mbox{, }D\less E\mbox{ is separated, and }
\forall F\in [E]^{<|E|}\;((D\less E)\cup F \mbox{ is separated})\}.
$$
Also in 
\cite{LL94}, it is shown that $\lambda^\omega<\theta$ implies that 
there are no 
$(\theta,\lambda)$-good sets; in particular,  $GCH$ implies that 
there are no $(\theta,\lambda)$-good sets whenever 
$\theta>\lambda\ge 
cf(\lambda)>\omega$.

In this paper, we give several forcing constructions of 
$(\theta,\lambda)$-good sets. Starting with a regular cardinal 
$\theta$, the first construction gives a model with a 
$(\theta,\omega)$-good set but no $(\theta,\omega)$-gaps (this 
model has been obtained independently by Haim Judah~\cite{Judah}). 
We then 
modify this construction to give a $(\theta,\omega)$-good set when 
$\theta$ is a singular cardinal of uncountable cofinality; this yields 
an easy consistent example of a first countable space in which cwH 
fails for the first time at a singular cardinal (different examples of 
such spaces were constructed by Fleissner and 
Shelah~\cite{Fleissner-Shelah} and by 
Koszmider~\cite{Koszmidor-Kurepa}). 
The final construction gives models with 
$(\theta,\lambda)$-good sets when $cf(\theta)\ge\omega_1$ and 
either $\lambda$ is regular and $\lambda\le \theta$ or 
$\lambda^\omega\le\theta$. When 
$\theta>\lambda$, the first countable $<\!\theta$-cwH not 
$\le\!\theta$-cwH spaces obtained from these good sets are new, in 
the sense that they can be made cwH by removing a small 
(cardinality $\lambda$) closed discrete set.

We also prove two theorems that characterize the existence of 
certain good sets in terms of integer-valued functions. These 
theorems will use the following relations $\le^+$ and $\le^*$ on 
$\omega^\lambda$ that 
generalize the usual notion of $\le^*$ on $\omega^\omega$.

\begin{Definition} Let $\lambda$ be an infinite cardinal, and let $f,g 
\in \omega^\lambda$. 
\begin{enumerate}
\item[(a)] We say $f \le^* g$ if for all but finitely many $\alpha <
\lambda$, $f(\alpha) \leq g(\alpha)$. 

\item[(b)] We say $f\le^+g$ if there is a $k\in \omega$ 
such that for all $\alpha<\lambda$, either $f(\alpha)\le g(\alpha)$ 
or $f(\alpha) \le k$.
\end{enumerate}
\end{Definition}

Notice that $\le^+$ is a reflexive, transitive relation on 
$\omega^\lambda$.  We put $f=^+ g$ if there is a $k\in \omega$ 
such that for all $\alpha\in \lambda$, either $f(\alpha) = g(\alpha)$ 
or both $g(\alpha),f(\alpha)\le k$. This 
determines an equivalence relation on $\omega^\lambda$, and   the 
order 
on these equivalence classes induced by $\le^+$ is a partial 
order.  If $\lambda = \omega$ and we restrict ourselves to strictly 
increasing functions, then the two notions $\le^*$ and 
$\le^+$ 
coincide.  Also note that $f\le^*g$ always implies $f\le^+g$.

The following lemma gives a canonical way to construct a 
$(\theta,\lambda)$-good set. We will prove two partial converses to 
this result later. Whenever we refer to a family of functions from 
some set $A$ into $\omega$ as being bounded or unbounded, we 
mean 
with respect to the obvious $\le^+$ order on $\omega^A$
(unless we state explicitly $\le^*$). When $\lambda =\omega$, the
two notions of being bounded or unbounded coincide.

\begin{Lemma}\label{functions-imply-good} Let $\lambda\le\theta$ 
be infinite cardinals. Assume 
there is $\Cal F = \{ f_\beta : \beta < \theta \}
\subseteq\omega^\lambda$ so that:
\begin{enumerate}
\item[\em{(a)}] $\Cal F$ is unbounded;
\item[\em{(b)}] each $\Cal  G\in [\Cal F]^{< \theta}$ is bounded;
\item[\em{(c)}] for all $E \in [\lambda]^{< \lambda}$,
$\Cal F \re E = \{ f_\beta \re E :\beta < \theta \}$
is bounded. 
\end{enumerate}
 Then there is a $(\theta,\lambda)$-good set. 
\end{Lemma}

\begin{Proof} Put $A = \{ (\langle \beta , m \rangle , \langle
\alpha , n \rangle ) : m,n \le f_\beta (\alpha) \}$. We show:
\begin{enumerate}
\item[(I)] $(\infty , \infty) \in \cl{A}$;

\item[(II)] $(\infty , \infty) \notin \cl{B}$ for all $B \in [A]^{
<\theta}$;

\item[(III)] $(\infty , \infty) \notin \cl{C}$ whenever $C =
A \cap ((\theta \times \omega) \times (E \times \omega))$
and $E \in [\lambda]^{<\lambda}$.
\end{enumerate}

\noindent (I) Choose $g \in \omega^\theta$ and 
$f\in\omega^\lambda$.
We have to show that $A$ intersects $V_g \times U_f$.
By (a), we can choose $\beta < \theta$ so that $f_\beta \not\leq^+
f$.  I.e.,  for every $k\in \omega$, there is an $\alpha_k<\lambda$ 
such that $f_\beta(\alpha_k) >f(\alpha_k)$ and 
$f_\beta(\alpha_k)>k$.  
Choose $k$ so that $ k \ge g (\beta)$.
Then $(\langle\beta , f_\beta (\alpha_k) \rangle ,\langle
\alpha , f_\beta (\alpha_k) \rangle ) \in A \cap (V_g \times U_f)$.

\smallskip
\noindent(II) Let $B \in [A]^{<\theta}$; let $G = \{ \beta <\theta:
\exists m, n <\omega\;\exists \alpha<\lambda \; 
(( \langle\beta, m 
\rangle , \langle \alpha , n \rangle ) \in B )\}$. Then 
$G \in [\theta]^{< \theta}$. By assumption (b),
there is an $f\in \omega^\lambda$ that is a $\le^+$-bound for $\Cal 
G = \{ f_\beta : \beta \in G \}$. Thus for all $  \beta \in 
G $, there is a 
$ k_\beta\in\omega$ such that for each $\alpha<\lambda$, either
$f_\beta (\alpha) \leq k_\beta$ or $f_\beta (\alpha) \leq 
f(\alpha)$.
Define $g\from\theta\to \omega$ by  $g(\beta) = k_\beta$ for 
$\beta \in G$, $g(\beta) = 0$ otherwise.
We have to show that $B \cap (V_g \times U_f) = \emptyset$.
To see this take $(\langle \beta , m \rangle , \langle \alpha , 
n\rangle ) \in
B$. Then both $m$ and $n$ are less than or equal to 
$f_\beta(\alpha)$; thus either $n \leq
f_\beta (\alpha) \leq f(\alpha)$ and $\langle \alpha , n \rangle
\notin U_f$ or $m \leq f_\beta(\alpha) \leq k_\beta = g(\beta)$
and $\langle \beta, m \rangle \notin V_g$.

\smallskip
\noindent(III) Similar; using (c) in place of (b).\end{Proof}

\section{ $(\theta,\omega)$-good sets, $(\theta,\omega)$-gaps, and 
incompactness at singulars, revisited\label{omega-good}}

Before we give a characterization of $(\theta,\omega)$-good sets in 
terms of families of integer-valued functions, we need to recall 
some facts about $(\theta,\lambda)$-good sets from \cite{LL94}. 
Given a  set $A \subseteq 
(\theta\times\omega)\times(\lambda\times \omega)$, and ordinals 
$\beta<\theta$ and $\alpha<\lambda$, we define 
$H_{\beta\alpha}(A) = H_{\beta\alpha} = 
\{(m,n):m,n\in\omega\mbox{ and 
}(\langle\beta,m\rangle,\langle\alpha,n\rangle
)\in A\}$. We say $H_{\beta\alpha}$ is {\em closed downward} 
(abbreviated cdw) if whenever $(m,n)\in H_{\beta\alpha}$, $n'\le n$, 
and $m'\le m$, then $(m',n')\in H_{\beta\alpha}$. The proofs of the 
following lemmas can be found in \cite{LL94}.

\begin{Lemma} \label{cdw} Suppose there is a 
$(\theta,\lambda)$-good set. Then 
there is a $(\theta,\lambda)$-good set $A$ such that each 
$H_{\beta\alpha}$ is finite and cdw.
\end{Lemma}

\begin{Lemma} \label{accumulation} Suppose $A\subseteq 
(\theta\times \omega)\times 
(\lambda\times \omega)$ and that each $H_{\beta\alpha}$ is finite 
and cdw. Then $A$ accumulates at $(\infty,\infty)$ in 
$F_\theta\times F_\lambda$ if and only if
$$
\forall f\from \lambda\to 
\omega\quad\exists\beta<\theta\quad\forall 
m\in\omega\quad\exists\alpha<\lambda\quad((\langle\beta,m
\rangle,\langle 
\alpha,f(\alpha)\rangle)\in A).
$$
\end{Lemma}

In \cite{G}, Gruenhage showed  that if there is a 
$(\theta,\omega)$-good 
set,  then $\theta\ge \frak b$.

\begin{Theorem} \label{char-omega-good} The following are 
equivalent:
\begin{enumerate}
\item[\em{(a)}] There is an $\Cal F\subseteq\omega^\omega$ of 
cardinality $\theta$ such that $\Cal F$ is $\le^*$-unbounded, but 
every $\Cal G\in [\Cal F]^{<\theta}$ is $\le^*$-bounded.
\item[\em{(b)}] There is a $(\theta,\omega)$-good set.
\end{enumerate}
\end{Theorem}
\begin{Proof} (a) $\Rightarrow$ (b) 
By Lemma~\ref{functions-imply-good}.

\noindent (b) $\Rightarrow$ (a) Suppose that $A$ is a 
$(\theta,\omega)$-good set. Without loss of generality, we can 
assume that each $H_{\beta k} = H_{\beta k}(A)$ is finite and cdw. 
By Gruenhage's result,  $\frak b\le \theta$; so let 
$\{g_\alpha:\alpha<\frak b\}\subseteq\omega^\omega$ be an 
unbounded family of strictly increasing functions.  For 
$\beta<\theta$ and  $\alpha<\frak b$, define a function 
$f_{\beta,\alpha}\from \omega\to \omega$ by
$$
f_{\beta,\alpha}(k) = \left\{ 
\begin{array}{ll}
\max\{n: \exists m\; (g_\alpha(m)\ge k \land 
(\langle\beta,m\rangle,\langle k,n\rangle)\in A)\}& {} \\
0 \quad\mbox{ if the above set is empty.}&{}
\end{array} \right.
$$

Set $\Cal F = \{f_{\beta,\alpha}: \beta<\theta $ and $\alpha<\frak 
b\}$; we check that (I) $\Cal F$ is $<^*$-unbounded and (II) every 
$\Cal G\in[\Cal F]^{<\theta}$ is $<^*$-bounded.

\medskip
\noindent (I) Let $f\in \omega^\omega$. Using 
Lemma~\ref{accumulation} and the fact 
that $A$ is $(\theta,\omega)$-good, there is a $\beta<\theta$ such 
that for each $m\in \omega$, there is a $k_m\in\omega$ such that 
$(\langle\beta,m\rangle,\langle k_m,f(k_m)\rangle)\in A$. Define 
$g\in \omega^\omega$ by $g(m) = k_m$; then $g$ is finite-to-one.
Because the $g_\alpha$'s are $\le^*$-unbounded, there is an 
$\alpha<\frak b$ such that for infinitely many $m$, $g_\alpha(m)\ge 
g(m)$.  Then for each such $m$, $f_{\beta,\alpha}(k_m)\ge f(k_m)$, 
so that $f$ is not a $\le^*$-bound for $\Cal F$.

\medskip
\noindent (II) Let $\Cal G\in[\Cal F]^{<\theta}$. Set $G = 
\{\beta<\theta: \exists \alpha<\frak b\;(f_{\beta,\alpha}\in \Cal 
G)\}$. Set $B = A\cap (G\times 
\omega)\times(\omega\times\omega)$; then $(\infty,\infty)$ is not 
in the closure of $B$. Choose $g\in \omega^\theta$ and $f\in 
\omega^\omega$ so that $B \cap (V_g\times U_f)= \emptyset$. We 
claim that $f$ is a $\le^*$-bound for $\Cal  G$.
Fix an $f_{\beta,\alpha} \in \Cal G$. Now, $\beta\in G$, so
$(\langle\beta,m\rangle,\langle k,n\rangle)\notin A$ whenever $m > 
g(\beta)$ and $n >f(k)$.  Stated contrapositively, if 
$(\langle\beta,m\rangle,\langle k,n\rangle)$ is in $A$, then either 
$m \le g(\beta)$ or $n \le f(k)$.  
Take $k>g_\alpha(g(\beta))$ and any $m$ such that $g_\alpha(m)\ge 
k$, then $g_\alpha(m)>g_\alpha(g(\beta))$.  Because $g_\alpha$ is 
strictly increasing, $m > g(\beta)$.  Thus, if $n$ is such that 
$(\langle\beta,m\rangle,\langle k,n\rangle)\in A$, we must have 
$n\le f(k)$. Taking the maximum over all such $n$ gives 
$f_{\beta,\alpha}(k)\le f(k)$. Thus, whenever $k> 
g_\alpha(g(\beta))$, we have $f_{\beta,\alpha}(k)\le f(k)$, whence 
$f_{\beta,\alpha}\le^* f$.\end{Proof}

It is a well-known fact, due independently to 
Hausdorff~\cite{Hausdorff} and Rothberger~\cite{Rothberger}, that 
the existence of a $(\theta,\omega)$-gap in 
$(\omega^\omega,\le^*)$ 
is equivalent to the existence of a well-ordered unbounded sequence 
of order type $\theta$ in $(\omega^\omega,\le^*)$.  This means that 
the 
existence of a $(\theta,\omega)$-gap implies the existence of a 
$(\theta,\omega)$-good set. It is certainly consistent that the 
converse is true, (consider, e.g., a model of $CH$), so it is natural to 
ask if the converse is true in $ZFC$.  In fact, when 
$\theta\ge\omega_2$ and $cf(\theta)\ge \omega_1$, there is a 
model in 
which there is a $(\theta,\omega)$-good set, but only  
$(\omega_1,\omega_1)$-  and 
$(\omega_1,\omega)$-gaps.

Before we define the partial orders that give these models, we state 
some 
useful lemmas about product forcing. If $\Cal F\subseteq 
\omega^\omega$ and $g\in \omega^\omega$, we say that $g>^*\Cal 
F$  if $g>^* f$ for all $f\in \Cal F$. We say that a real $f$ in a 
universe 
larger than $V$ is {\em unbounded over $V$\/} if $f\not\le^*g$ for 
all  $ g\in\omega^\omega\cap V$.

\begin{Lemma}\label{unbounded} Let $\bbb P$ and $\bbb Q$ be 
partial orders. Suppose 
$\dot f $ is a $\bbb P$-name for a real and $\dot g$ is a $\bbb Q$-
name for a real.  If $\forces_{\bbb P}``\dot f$ is unbounded over 
$V$'', then $\forces_{\bbb P\times \bbb Q}``\neg(\dot f\le^* \dot 
g)$''.
\end{Lemma}

\begin{Proof} Fix $n\in\omega$ and $(p,q)\in \bbb P\times \bbb Q$. 
Because $\forces_{\bbb P}``\dot f$ is unbounded over $V$'', there is 
an $l\ge n$ such that $\forall m\in \omega$, $\neg(p\forces_{\bbb 
P}``\dot f(l)\le m")$.

Choose a $q'\le q$ and an $m\in \omega$ so that $q'\forces_{\bbb 
Q}``\dot g(l) = m$''. Then we can find a $p'\le p$ and a $k> m$ such 
that $p' \forces_{\bbb P}``\dot f(l) = k$''. Thus, $(p',q')\forces_{\bbb 
P\times \bbb Q}``\dot f(l)>\dot g(l)$''.\end{Proof}

\begin{Lemma} \label{no-interpolation} Suppose that $\dot h_1$ and 
$\dot h_2$ are $\bbb 
P$-names for reals, $\dot g$ is a $\bbb Q$-name for a real, and
$(p,q)\forces_{\bbb P\times \bbb Q} ``\dot h_1<^* \dot g <^* \dot 
h_2$''. Then $p\forces_{\bbb P}``\exists j\in \omega^\omega\cap 
V\;(\dot h_1<^* j <^* \dot h_2)$''.\end{Lemma}

\begin{Proof} Without loss of generality, we assume that for some 
$k'\in \omega$, 
$$(p,q)\forces_{\bbb P\times \bbb Q}``\forall n\ge 
k'\;(\dot h_1(n)<^* \dot g(n) <^* \dot h_2(n))".$$

For each $n\in\omega$, define $j(n)$ by taking a $q'\le q$ and an 
$m\in \omega$ so that $q'\forces_{\bbb Q}``\dot g(n) = m$'', and 
setting $j(n) = m$.

We claim that $p\forces_{\bbb P} ``\dot h_1<^* j <^* \dot h_2$''. 
Otherwise, for each $k\in \omega$, there is an $n\ge k$ and a $p'\le 
p$ such that $p'\forces_{\bbb P}``\neg(\dot h_1(n)< j(n) < \dot 
h_2(n))$''. Fix such a $p'$ and $n$ for the $k'$ given above, and find a 
$q'\le q $ such that 
$q'\forces_{\bbb Q}``\dot g(n) = j(n)$''. Then
$$
(p',q')\forces_{\bbb P\times \bbb Q}``(\dot h_1(n)< \dot g(n) < \dot 
h_2(n)) \land \neg (\dot h_1(n)<   j(n) < \dot h_2(n))\land  (\dot 
g(n) = j(n))",$$
a contradiction.\end{Proof}

For a set $A$, $\bbb C_A$ is the partial order $\Fn(A\times 
\omega,\omega) = \{p:p$ is a finite partial function from $A\times 
\omega$ into $\omega\}$. $V_A$ is the generic extension of $V$ by 
$\bbb C_A$. We will need the fact, due to Kunen~\cite{Kunen-diss}, 
that forcing with $\bbb C_A$ over a model of $CH$ does not add an 
$\omega_2$-sequence in $(\omega^\omega,\le^*)$.

\begin{Lemma}\label{no-chains} Assume $CH$ and set $\theta = 
\omega_2$. In 
$V_\theta$, let $\bbb P = \prod_{\alpha<\theta} \bbb P_\alpha$ be a 
ccc finite support product of $\aleph_1$-sized partial orders such 
that $\prod_{\alpha<\beta} \bbb P_\alpha\in V_\beta$. Set $\bbb Q = 
\bbb C_\theta\star\dot\bbb P$, and let $H$ be $\bbb Q$-generic over 
$V$. 
Then there are no well-ordered $\omega_2$-sequences in 
$V[H]$.\end{Lemma}

\begin{Proof} By way of contradiction, let $\{\dot 
f_\alpha:\alpha<\omega_2\}$ be a collection of $\bbb Q$-names for 
a well-ordered $\omega_2$-sequence. Without loss of generality, 
each $\dot f_\alpha = \bigcup_{n,m\in\omega}\{(m,n)\}\times 
A_{mn}$, where $A_{mn}$ is a maximal countable antichain and each 
$q\in A_{mn}$ has the form $q = (c, \langle \dot p_{\alpha_1}, \dot 
p_{\alpha_2},\dots \dot p_{\alpha_k}\rangle)$, where $c \in 
\bbb C_\theta$ and $\dot p_{\alpha_i}\in \dot \bbb P_{\alpha_i}$. 
I.e., because 
the supports of conditions in $\dot\bbb P$ are finite, we can assume 
that the $\bbb C_\theta$ part of a condition is strong enough to 
decide 
the support of the $\dot\bbb P$ part.

Define $\suppt(q) =  \{\alpha_1, 
\alpha_2,\dots,\alpha_k\}$ and for each $\alpha\in \omega_2$, 
$$A_\alpha = \suppt(\dot f_\alpha) = 
\bigcup_{\stackrel{m,n\in\omega}{q\in 
A_{mn}}}\suppt(q).$$

By   thinning and re-indexing, we can assume that the 
$A_\alpha$'s are a delta system with root $\Delta$ and that 
$\alpha<\alpha'<\omega_2$ implies
$$
\max \Delta < \min (A_\alpha\less \Delta) \le \max (A_\alpha\less 
\Delta) < \min (A_\alpha'\less \Delta).
$$
Set $\beta = \max(\Delta)+1$, and force with 
$\bbb C_\beta\star\prod_{\alpha<\beta} \dot\bbb P_\alpha$. Notice 
that 
$CH$ 
is still true, so if we force with $\bbb C_{[\beta,\theta)}$, we 
obtain a 
model $V'$ with no $\omega_2$-chains.  We can also assume that 
each $\dot f_\alpha$ is a 
$\prod_{\xi\in A_\alpha\less \Delta}\bbb P_\xi$-name.

Set $E = \{\xi\in \omega_2\less \beta:$ For some even $\alpha$, 
$\xi\in A_\alpha\less\Delta$.$\}$, and set $O = 
\omega_2\less(E\cup 
\beta)$. Finally, let $\bbb P_E = \prod_{\xi\in E} \bbb P_\xi$ and
$\bbb P_O = \prod_{\xi\in O}\bbb P_\xi$.

Working in $V'$, fix a condition $(p,q) \in \bbb P_E\times \bbb P_O$ 
such that 
$$
(p,q) \forces_{\bbb P_E\times \bbb P_O}`` \forall 
\alpha<\omega_2\;(\alpha\mbox{ even implies }
\dot f_\alpha<^* \dot f_{\alpha+1}<^*\dot f_{\alpha+2})".
$$
By the Lemma~\ref{no-interpolation}, we can find for each even 
$\alpha<\omega_2$, a $p_\alpha\le p$ and a 
$j_\alpha\in \omega^\omega\cap V'$ such that 
$$
p_\alpha\forces_{\bbb P}``\dot f_\alpha<^* j_\alpha<^*\dot 
f_{\alpha+2}".
$$
Now, because supports are finite and each $|\bbb P_\alpha|\le 
\omega_1$, we can find an $A\in [\omega_2]^{\omega_2}$ so 
that whenever $\alpha,\alpha'\in A$, then $p_\alpha$ and 
$p_{\alpha'}$ 
are 
compatible. But then $\alpha<\alpha'\in A$ implies that 
$j_\alpha<^*j_{\alpha'}$, contradicting the fact that there are no 
$\omega_2$-sequences in $V'$.\end{Proof}

 Hechler forcing is the partial order $\bbb D = \{(s,f): 
s\in \omega^{<\omega}$, $f\in\omega^\omega$, and $s\subseteq 
f\}$, ordered so that $(s,f)\le (t,g)$ if and  only if $ s\supseteq t$ 
and $\forall n\in\omega\;(f(n)\ge g(n))$.  Clearly, $\bbb D$ is 
$\sigma$-centered and adds a real that eventually dominates all 
ground model reals. The following result has been obtained 
independently by Judah~\cite{Judah}.

\begin{Theorem} Assume $V\models CH$ and  set $\theta= 
\omega_2$.  Then there is a partial 
order $\bbb Q$ such that whenever $H$ is $\bbb Q$-generic over $V$, 
then the following are true in $V[H]$:
\begin{enumerate}
\item[{\em (1)}] There is a family $\Cal F = 
\{f_\alpha:\alpha<\theta\}\subseteq 
\omega^\omega$ such that:
\begin{enumerate}
\item[{\em (a)}] for all $\beta<\theta$, $\{f_\alpha:\alpha<\beta\}$ 
is bounded;
\item[{\em (b)}] $\Cal F$ is unbounded.
\end{enumerate}
\item[{\em (2)}] There are no well-ordered sequences of length 
$\omega_2$ in 
$(\omega^\omega,\le^*)$.
\end{enumerate}
\end{Theorem}

\begin{Proof} We define $\bbb Q$ as a two-step iteration. The first 
step of the iteration is simply $\bbb C_\theta$.  We define the 
second step by working in $V_\theta$.  Let $\bbb P_\beta$ be $\bbb 
D^{V_\beta}$, i.e., Hechler forcing in the sense of the model obtained 
by adding the first $\beta$-many Cohen reals.  In $V_\theta$, each 
$\bbb D_\beta$ is $\sigma$-centered, so the finite support product 
$\bbb P = \prod_{\beta<\theta}\bbb P_\beta$ is ccc. Hence, $\bbb Q 
= 
\bbb C_\theta\star \bbb P$ is ccc.

Let $H$ be $\bbb Q$-generic; in $V[H]$, define $\Cal F = 
\{f_\beta:\beta<\theta\}$, where $f_\beta$ is the $\beta$-th 
Cohen real.  Let $g_\beta$ be the Hechler real added by $\bbb 
P_\beta$ over $V_\theta$; because 
$\{f_\alpha:\alpha<\beta\}\subseteq V_\beta$, we have 
$f_\alpha\le^* g_\beta$ for each $\alpha<\beta$. This establishes 
(a).

Let $\dot g$ be a $\bbb Q$-name for a real. Note that $\forall 
\beta<\theta$, we have 
$$
\bbb Q \cong [(\bbb C_\beta\star\prod_{\alpha<\beta}\dot\bbb 
P_\alpha)\times \bbb C_{[\beta,\theta)}]\star \prod_{\alpha\ge 
\beta}\dot \bbb P_\alpha.
$$
By this observation and the fact that  $\bbb Q$ is ccc, there is a 
$\beta<\theta$ such that $\dot g$ is a $\bbb C_\beta\star 
\prod_{\alpha<\beta}\bbb P_\alpha$-name.  Because $f_\beta$ is 
unbounded over $V$, Lemma~\ref{unbounded} and the 
fact that $\le^*$ is upwards absolute imply that $\forces_{\bbb Q}$ 
``$\dot f_\beta \not\le^* \dot g$''. Hence, no real in $V[H]$ bounds 
$\Cal F$, so (b) is established.
 
Note that the partial order we've defined satisfies the hypotheses of 
Lemma~\ref{no-chains},  so (2) is 
also true.\end{Proof}

We can also show, via a modification of the  ``isomorphism of 
names'' argument originally due to Kunen~\cite{Kunen-diss}, that the 
forcing construction given above yields a model with no 
$\omega_2$-sequences when $\theta$ is any cardinal (see \cite{Br}
for details).
When $\theta$ is regular, this gives a model with a 
$(\theta,\omega)$-good set, but only $(\omega_1,\omega_1)$- and 
$(\omega_1,\omega)$-gaps.

The above proof does not quite suffice to produce a 
$(\theta,\omega)$-good set when $\theta$ is a singular 
cardinal of uncountable cofinality, because we need to 
bound all small subfamilies of $\Cal F$.  Fortunately,  we can use 
the ccc to accomplish this.

\begin{Theorem} Suppose $V\models GCH$ and that $\omega_1\le 
cf(\theta)<\theta$.  Then there is a ccc partial order $\bbb Q$ such 
that whenever $H$ is $\bbb Q$-generic over $V$, there is a family 
$\Cal F = \{f_\alpha:\alpha<\theta\}\subseteq\omega^\theta$ in 
$V[H]$  such that:
\begin{enumerate}
\item[{\em (1)}] $\Cal F$ is unbounded;
\item[{\em (2)}] $\forall \Cal G \in [\Cal F]^{<\theta}$, $\Cal G$ is 
bounded.
\end{enumerate}
\end{Theorem}
\begin{Proof} As before, $\bbb Q$ will be a two step iteration, with 
first step $\bbb C_\theta$. Let 
$f_\beta$ be the $\beta$-th Cohen real, and set $\Cal F = 
\{f_\beta:\beta<\theta\}$.  The second step of the iteration is 
defined in $V_\theta$ as $\prod_{A\in 
[\theta]^{<\theta}}\bbb P_A$, where $\bbb P_A = \bbb D^{V_A}$, i.e., 
Hechler forcing in the sense of the model obtained by adding the 
Cohen reals with indices in $A$.   We therefore have
$$\bbb Q = \bbb C_\theta\star \prod_{A\in 
[\theta]^{<\theta}}\dot \bbb P_A.$$
  Clearly, $\bbb Q$ is ccc.  As before, $\Cal F$ remains unbounded in 
$V[H]$.

To see that every small family is bounded, take $B\in 
[\theta]^{<\theta}\cap V[H]$.  Because $\bbb Q$ is ccc, there is an 
$A\in [\theta]^{<\theta}\cap V$ such that $B\subseteq A$.  Then the 
Hechler real added by $\bbb P_A$ bounds $\{f_\beta:\beta\in 
B\}$.\end{Proof}

Again it can be shown that the model does not contain
well--ordered sequences of length $\omega_2$ in $(\omega^\omega ,
\le^*)$.

Via the translation results in \cite{LL94}, this theorem 
gives a new example of a first countable space in which cwH fails 
for the first time at a singular cardinal. Notice that this space can 
be made 
cwH by removing a countable closed discrete set---a property that 
previous examples of such spaces did not have.

\section{Forcing $(\theta,\lambda)$-good sets}

In this section, we describe notions of forcing for adding 
$(\theta,\lambda)$-good sets for some cardinals that satisfy 
$\omega_1\le \lambda\le\theta$.  The method will be similar to 
that used in the previous section---but proving that our iteration is 
ccc will now be non-trivial.

\begin{Theorem} \label{main-result} Let $\lambda\le \theta$ 
cardinals with $cf(\theta)
\ge \omega_1$. Assume either $\lambda$ is regular or 
$\lambda^\omega
\le \theta$. There
is a ccc partial order $\bbb P$ that adds a family $\Cal 
F=\{f_\xi:\xi<\theta\}\subseteq 
\omega^\lambda$ satisfying:
\begin{enumerate}
\item[\em{(a)}] $\Cal F$ is $\le^+$--unbounded;
\item[\em{(b)}] for all $ B \subseteq \theta$ with $|B| 
<\theta$, $\{f_\xi:\xi\in B\}$ is 
$\le^*$--bounded;
\item[\em{(c)}] for all $A\subseteq \lambda$ with $|A| <\lambda$,
$\Cal F \re A = \{ f_\xi \re A :\xi < \theta \}$
is $\le^*$--bounded. 
\end{enumerate}\end{Theorem}
\begin{Proof} Let $V_0$ be the ground model.  As before, $\bbb P\in 
V_0$ will be a two-step iteration.  To define the first step, we 
start with a function $H\from \theta\to[\lambda]^\omega$ that 
satisfies either
\begin{enumerate}
\item[(1)] $\forall A\in [\lambda]^\omega \; (|\xi <\theta :
H(\xi) = A \} | = \theta)$
\end{enumerate}
or (in case $\lambda$ is regular and $\lambda^\omega \leq \theta$
fails)
\begin{enumerate}
\item[(2)] $\forall \xi<\theta\;\forall \zeta<\lambda\;(\zeta\in 
H(\xi)\Rightarrow \zeta+1\in H(\xi))$ and
\item[(3)] $\forall \zeta<\lambda\;(cf(\zeta) = \omega\Rightarrow 
|\{\xi<\theta:\sup(H(\xi))=\zeta\}| = \theta)$.
\end{enumerate}

We define $\bbb P_0 = \{c\in \bbb C_{\theta\times\lambda}: \forall 
\xi<\theta\;\forall \zeta<\lambda\;((\xi,\zeta)\in \dom(c) 
\Rightarrow \zeta\in H(\xi))\}$, ordered by reverse containment
(notice that our notation here is slightly different from section
2: $\bbb C_A$ denotes $Fn (A,\omega)$). 
Obviously, $\bbb P_0$ is forcing isomorphic to $\bbb C_\theta$; we 
think 
of $h_\xi$ (the $\xi$-th Cohen real added by $\bbb P_0$) as 
having domain $H(\xi)$.  Extend $h_\xi$ to a function 
$f_\xi\in \omega^\lambda$ by setting $f_\xi(\zeta) = h_\xi(\zeta)$ 
for 
all $\zeta\in H(\xi)$ and $f_\xi(\zeta) = 0$ for all 
$\zeta\notin H(\xi)$.

Let $V_1$ be the extension of $V_0$ by $\bbb P_0$. In $V_1$, we 
define for each $A\subseteq \lambda$ with $|A|<\lambda$ and $A \in 
V_0$
and for each $B\subseteq\theta$ with $|B| <\theta$ and $B\in V_0$
partial orders 
$\bbb Q_A$ and $\bbb R_B$ as follows:  
$\bbb Q_A = \{\langle s,F\rangle: s\in \bbb C_A\land F\in 
[\theta]^{<\omega}\}$ ordered so   that $\langle s,F\rangle\le 
\langle s',F'\rangle$ if and  only if $s \supseteq s'$, $F \supseteq
F'$  and 
$$
\forall 
\xi\in F'\;\forall \zeta\in\dom(s)\less\dom(s')\;( 
f_\xi(\zeta)\le s(\zeta) ).
$$
Similarly, $\bbb R_B =\{\langle t,G\rangle:t \in \bbb C_\lambda 
\land G\in [B]^{<\omega}\}$, ordered in the same way as $\bbb Q_A$: 
$\langle 
t,G\rangle\le \langle t',G'\rangle$ if and only if $t\supseteq t'$,
$G \supseteq G'$  and 
$$
\forall 
\xi\in G'\;\forall \zeta\in\dom(t)\less\dom(t')\;( 
f_\xi(\zeta)\le t(\zeta) 
).
$$
 
In $V_1$, let $\bbb P_1 = 
\prod_A\bbb Q_A \times 
\prod_B\bbb R_B$ be the finite support product of 
the 
$\bbb Q_A$'s and $\bbb R_B$'s. In $V_0$, let $\bbb P = \bbb 
P_0\star\dot\bbb P_1$. Also in $V_0$, let $\langle A_\alpha : 
\alpha < \lambda '\rangle$ enumerate 
$[\lambda]^{<\lambda}$ and $\langle B_\beta : \beta < \theta ' 
\rangle
$ enumerate $[\theta]^{<\theta}$.

\medskip
\noindent{\bf Claim:\/} $\bbb P_1$ is ccc in $V_1$.

To prove the claim, it suffices to show that for every 
$A\in[\lambda ']^{<\omega}$ and $B\in [\theta ']^{<\omega}$, the 
partial 
order
$\bbb P_0 \star (\prod_{\alpha\in A}\dot \bbb 
Q_{A_\alpha}\times \prod_{\beta\in B}\dot \bbb R_{B_\beta})$ is 
ccc in 
$V_0$.  
For each $\gamma\in \omega_1$, fix a condition
$$
p^\gamma = \langle 
c^\gamma,\langle\langle s^\gamma_\alpha,F^\gamma_\alpha\rangle: 
\alpha\in A\rangle,\langle\langle t^\gamma_\beta, 
G^\gamma_\beta\rangle:\beta\in B\rangle\rangle.
$$
We don't need to work with names 
because conditions are finite partial functions, so we can assume 
the $\bbb P_0$ part of a condition is strong enough to decide the 
second part. For $c\in \bbb P_0$, set $d(c) = \{\xi<\theta:\exists 
\zeta\in 
H(\xi)\;((\xi,\zeta)\in\dom(c))\}$. 

By applying a delta-system argument, we can assume that there are 
$c$, $s_\alpha$, $F_\alpha$, $t_\beta$, and $G_\beta$ such that 
$\forall \gamma<\omega_1$, $\forall \alpha\in A$, and $\forall 
\beta\in B$,
$$p^\gamma = \langle c\cup 
c^\gamma,\langle\langle s_\alpha\cup s^\gamma_\alpha,F_\alpha
\cup F^\gamma_\alpha\rangle:\alpha\in 
A\rangle,\langle\langle t_\beta\cup t^\gamma_\beta,G_\beta\cup 
G^\gamma_\beta\rangle:\beta\in B\rangle\rangle$$
 and $\forall \gamma<\delta<  \omega_1$, $\forall \alpha\in 
A$, and $\forall \beta\in B$, each of $d(c^\gamma)\cap d(c^\delta) $,  
$d(c) \cap d(c^\gamma)$,
$\dom(s^\gamma_\alpha)\cap\dom(s^\delta_\alpha)$, 
$F^\gamma_
\alpha\cap F^\delta_\alpha$, $ 
\dom(t^\gamma_\beta)\cap\dom(t^\delta_\beta)$, and $
G^\gamma_\beta\cap G^\delta_\beta$ is the empty set (this also
uses the countability of the sets $H(\xi)$).

We now define $\forall \gamma<\omega_1$,
\begin{eqnarray*}
P(\gamma) &= & d(c^\gamma)\cup\bigcup_{\alpha\in 
A}F^\gamma_\alpha\cup\bigcup_{\beta\in 
B}G^\gamma_\beta,\\
Q(\gamma) &=& \bigcup_{\alpha\in 
A}\dom(s^\gamma_\alpha)\cup\bigcup_{\beta\in 
B}\dom(t^\gamma_\beta)\mbox{, and }\\
P &= &d(c) \cup \bigcup_{\alpha\in 
A}F_\alpha\cup\bigcup_{\beta\in B}G_\beta.
\end{eqnarray*}

It is now easy to find $\gamma<\delta<\omega_1$ such that 
\begin{enumerate}
\item[(i)] $\emptyset =P(\gamma)\cap P = P\cap P(\delta) = 
P(\gamma)\cap P(\delta)$,
\item[(ii)] $Q(\gamma) \cap \{\zeta<\lambda:\exists \xi\in 
d(c^\delta)\cup d(c)\;((\xi,\zeta)\in \dom(c^\delta))\} = 
\emptyset$, and
\item[(iii)] $Q(\delta) \cap \{\zeta<\lambda:\exists \xi\in 
d(c^\gamma)\cup d(c)\;((\xi,\zeta)\in\dom(c^\gamma))\} = 
\emptyset$.
\end{enumerate}

We claim that $p^\gamma$ and $p^\delta$ are compatible. To see 
this, note that by (i), (ii), and (iii), we can find a $\hat c\supseteq 
c^\gamma\cup c^\delta$ such that
\begin{enumerate} 
\item [$(*)$] If $\xi\in P(\gamma)\cup P$ and $\zeta \in H(\xi)\cap 
Q(\delta)$, then $\langle\xi,\zeta\rangle\in \dom(\hat c)$ and 
$\hat c(\xi,\zeta) = 0 $.
\item[$(**)$]If $\xi\in P(\delta)\cup P$ and $\zeta\in H(\xi)\cap 
Q(\gamma)$, 
then $\langle\xi,\zeta\rangle\in \dom(\hat c)$ and $\hat 
c(\xi,\zeta) = 0$.
\end{enumerate}

Consider the condition 
$$p = \langle\hat c,\langle\langle s_\alpha\cup 
s^\gamma_\alpha\cup s^\delta_\alpha,F_\alpha\cup 
F^\gamma_\alpha
\cup F^\delta_\alpha\rangle:\alpha\in A\rangle,\langle\langle 
t_\beta\cup t^\gamma_\beta\cup t^\delta_\beta, G_\beta\cup 
G^\gamma_\beta\cup G^\delta_\beta\rangle:\beta\in 
B\rangle\rangle;$$
 we show that $p\le p^\gamma$ ($p\le p^\delta$ is 
similar).  

Clearly, all inclusion relations are met. Notice that by $(*)$ and 
$(**) $ we have  $\forall \alpha\in A$, $\forall \zeta\in\dom(
s^\delta_\alpha )$, and $\forall \xi\in F_\alpha\cup 
F^\gamma_\alpha$,
$$
0=\hat c(\xi,\zeta)\le 
s^\delta_\alpha(\zeta)\mbox{, or } 
 \zeta\notin H(\xi),
$$
and $\forall \beta\in B $, $\forall \zeta\in \dom(t^\gamma_\beta)$, 
and $\forall \xi\in G_\beta\cup 
G^\gamma_\beta$,
$$
0= \hat c(\xi,\zeta)\le t^\delta_\beta(\zeta)\mbox{, or }  
\zeta\notin H(\xi),
$$
so that $p\le p^\gamma$. This establishes the claim.

\medskip
Let $V_2$ be the extension of $V_1$ by $\bbb P_1$; notice that (b) 
and (c) of the theorem are true by genericity.  To complete the proof 
of the theorem, we need to show:

\medskip
\noindent{\bf Claim:} $\Cal F$ is unbounded in $V_2$.

By way of contradiction, suppose that there is a $\bbb P$-name 
$\dot f$ for an element of $\omega^\lambda$ such that
$\forces_{\bbb P}``\forall \xi<\theta\;(\dot f_\xi\le^+\dot f)$". (We 
will see later that $\bbb P$ factors nicely, so if this statement is 
only forced by some non-trivial condition $p$, we can replace $V_0$ 
with an initial extension obtained from a generic that contains $p$, 
and then argue as below.)

Assume first $cf (\lambda) \geq \omega_1$ and $\lambda^\omega
\leq \theta$.
Using condition (1) of the function $H$,  construct, by 
recursion on $\gamma<\omega_1$, conditions
$
p^\gamma = \langle 
c^\gamma,\langle\langle s^\gamma_\alpha,F^\gamma_\alpha\rangle: 
\alpha\in A^\gamma\rangle,\langle\langle t^\gamma_\beta, 
G^\gamma_\beta\rangle:\beta\in B^\gamma\rangle\rangle$
where
$$
c^\gamma \forces `` \langle s_\alpha^\gamma , F_\alpha^\gamma 
\rangle
\in \dot \bbb Q_{A_\alpha}  \mbox{ for } \alpha \in A^\gamma 
\mbox{
and } \langle t_\beta^\gamma , G_\beta^\gamma \rangle \in \dot 
\bbb R_{
B_\beta} \mbox{ for } \beta \in B^\gamma ", 
$$
ordinals $\alpha^\gamma<\lambda$ and $\beta^\gamma<\theta$, and 
integers $k^\gamma<\omega$ such that if 
$$
A(\gamma) = \bigcup_{\alpha \in A^\gamma} A_\alpha \mbox{ and }
B(\gamma) = \bigcup_{\beta \in B^\gamma} B_\beta,
$$
then
\begin{enumerate}
\item[(i)] if $cf(\gamma) = \omega$, then $\alpha^\gamma = 
\sup\{\alpha^\delta:\delta<\gamma\}$;
\item[(ii)] if $\gamma$ is a successor ordinal, then 
$\forall \delta <\gamma\; (\alpha^\gamma > \alpha^\delta
$ and $\alpha^\gamma \notin A(\delta))$;
\item[(iii)]   $\forall \delta<\gamma\;(\beta^\gamma  
\notin \{ \beta^\delta \} \cup B(\delta) \cup d(c^\delta) )$;
\item[(iv)] $H(\beta^\gamma)\supseteq
\{ \alpha^\delta : \delta < \gamma \}$; and
\item[(v)] $p^\gamma\forces_{\bbb P}``\forall \zeta<\lambda\;(\dot 
f_{\beta^\gamma}(\zeta)\le k^\gamma$ or $\dot 
f_{\beta^\gamma}(\zeta)\le \dot f(\zeta))$''.
\end{enumerate}

To avoid having to work with names, we are again assuming that
the $\bbb P_0$ part of a condition is strong enough to decide the 
$\bbb P_1$ part. Also notice that by conditions  (i) 
and (ii) of 
the recursion, 
$S_0= 
\{\alpha^\gamma:\gamma<\omega_1\}$ is club in $\alpha '
= \sup \{ \alpha^\gamma : \gamma <\omega_1\}$.
Thus $S \subseteq S_0$ is stationary in $\alpha'$ if and only if
$\tilde S = \{ \gamma<\omega_1 : \alpha^\gamma \in S\}$ is 
stationary
in $\omega_1$.

We now define regressive functions $a,b,c: \omega_1 \to 
[\omega_1]^{
<\omega}$ and $k :\omega_1 \to \omega$ by:
\begin{eqnarray*}
a(\gamma) &= &\{\delta<\gamma:
(A^\gamma \cap A^\delta) \less \bigcup_{\epsilon <\delta} 
A^\epsilon
\ne \emptyset \} ;\\
b(\gamma)&=& \{\delta<\gamma:
(B^\gamma \cap B^\delta) \less \bigcup_{\epsilon <\delta} 
B^\epsilon
\ne \emptyset \} ;\\
c(\gamma)&=& \{\delta<\gamma: (d(c^\gamma) \cap d(c^\delta) ) 
\less
\bigcup_{\epsilon <\delta} d(c^\epsilon) \ne \emptyset \} ;
\mbox{ and}\\
k(\gamma)&=&k^\gamma.
\end{eqnarray*}
By Fodor's lemma, we can find $\Delta_a , \Delta_b , \Delta_c \in
[\omega_1]^{<\omega}$, a 
$k_0\in \omega$, and a stationary $\tilde S_1\subseteq \omega_1$ 
such that  
$\forall  \gamma\in \tilde S_1$, $a(\gamma) = \Delta_a$, 
$b(\gamma) = \Delta_b$, $c(\gamma) = \Delta_c$ and 
$k(\gamma) = k_0$.

Set $\bar A = \bigcup_{\gamma \in \Delta_a} A(\gamma)$ and $\bar 
B =
\bigcup_{\gamma \in \Delta_b} B(\gamma) \cup \bigcup_{\gamma 
\in \Delta_c}
d(c^\gamma)$.
Then for $\gamma\neq \gamma' \in \tilde S_1$, we have
\begin{enumerate}
\item[(a)] $(A^\gamma\cap A^{\gamma'}) \less \bigcup_{\delta \in
\Delta_a} A^\delta = 
\emptyset$,
\item[(b)] $(B^\gamma\cap B^{\gamma'}) \less \bigcup_{\delta \in
\Delta_b} B^\delta = 
\emptyset$, and
\item[(c)] $(d(c^\gamma)\cap d(c^{\gamma'})) \less \bar B = 
\emptyset$.
\end{enumerate}

We now factor $\bbb P_0$ as $\bbb P_0^0\times \bbb P^1_0$, where
$$\bbb P^0_0=\{c\in\bbb P_0:\forall 
\zeta<\lambda\;\forall\xi<\theta\;((\xi,\zeta)\in 
\dom(c)\Rightarrow \xi\in \bar B \lor \zeta\in \bar A)\}
$$
and
$$
\bbb P^1_0 = \{c\in\bbb P_0:\forall 
\zeta<\lambda\;\forall\xi<\theta\;((\xi,\zeta)\in 
\dom(c)\Rightarrow \xi\notin \bar B \land \zeta\notin\bar A)\}.
$$
In turn, if we set $A = \bigcup_{\delta\in \Delta_a}A^\delta$, $B = 
\bigcup_{\delta\in \Delta_b}B^\delta$,
$$
\bbb P^0_1 = \prod_{\alpha\in A} \bbb Q_{A_\alpha}\times 
\prod_{\beta\in B} \bbb 
R_{B_\beta}\quad\mbox{and}\quad
\bbb P^1_1 = \prod_{\alpha\notin A}
\bbb Q_{A_\alpha}\times 
\prod_{\beta\notin B}\bbb R_{B_\beta},
$$
then $\bbb P$ can be factored as
$$
\bbb P = (\bbb P^0_0\star \dot \bbb P^0_1)\star (\bbb P^1_0 \star 
\dot \bbb P^1_1).
$$

Let $G^0$ be $\bbb P^0 = \bbb P^0_0\star \dot\bbb P^0_1$-generic 
over 
$V_0$, and let $V^0 = V_0[G^0]$. In $V^0$, set $\bbb P^1 = \bbb 
P^1_0\star\dot\bbb P^1_1$.

We claim that $G^0$ can be chosen so that $\tilde S^0 = 
\{\gamma\in \tilde S_1: p^\gamma\re\bbb P^0\in G^0\}$ is a 
stationary subset of $\omega_1$ in $V^0$. Otherwise, we can find a 
$\bbb P^0$-name $\dot C$ for a club subset of $\omega_1$ such that 
$\forces_{\bbb P^0}``\forall \gamma<\omega_1\;(\gamma\in 
\dot C\Rightarrow p^\gamma\re\bbb P^0\notin \dot G^0)$''. Because 
$\bbb 
P^0$ is ccc, we can find a club $C_0 \in V_0$ such that 
$\forces_{\bbb P^0}``C_0\subseteq\dot C$''.  Take a 
$\gamma\in C_0\cap \tilde S_1$, then
$$
p^\gamma\re\bbb P^0\forces_{\bbb P^0}``p^\gamma\re\bbb P^0\in 
\dot G^0\mbox{ and } \gamma\in \dot C",
$$
a contradiction.

Let $S^0 = \{\alpha^\gamma:\gamma\in \tilde S^0\}$. For 
$ \gamma\in \tilde S^0$, define 
$$
e(\alpha^\gamma) = \max \{ \zeta < \alpha^\gamma : (\beta^\gamma ,
\zeta) \in \dom (c^\gamma ) \}.
$$
(When we are talking about $c^\gamma$ (or other parts of 
conditions), we really mean 
$c^\gamma\re\bbb P^1_0$---this is ok because we've chosen 
$c^\gamma$ so that $c^\gamma\re \bbb P^0_0\in G^0$.) Notice that 
$e$ 
is regressive on $S^0$, so there is a stationary $T^0\subseteq S^0$ 
and a $\zeta_0<\alpha '$ so that $\forall \gamma\in 
T^0\;(e(\gamma) = \zeta_0)$. Set $\tilde T^0 = \{\gamma<\omega_1: 
\alpha^\gamma\in T^0\}$.

Now choose $\delta_0<\omega_1$ so that $\alpha^{\delta_0}
\notin \bar A$ and $\alpha^{\delta_0} > \zeta_0$.
Recall that $\alpha^{\delta_0} \in H(\beta^\gamma)$ whenever
$\gamma > \delta_0$ and $\gamma \in \tilde T^0$.
Also notice that without loss $\beta^\gamma \notin \bar B$ for
$\gamma \in \tilde T^0$.
Let $\dot G^1$ 
be the canonical name for a $\bbb P^1$-generic filter over $V^0$. We 
claim that
$$
(*)\qquad \forces_{\bbb P^1}``\forall k\in \omega\;\exists 
\gamma\in \tilde T^0\;(p^\gamma\in \dot G^1 \land \dot 
f_{\beta^\gamma}(\alpha^{\delta_0}) = k)".
$$
To see this, suppose that $p^1\in \bbb P^1$ and $k\in \omega$. Say
$$p^1 = \langle c^1,\langle\langle 
s^1_\alpha,F^1_\alpha\rangle:\alpha\in 
A^1\rangle,\langle\langle t^1_\beta,G^1_\beta\rangle:\beta\in 
B^1\rangle\rangle.$$

By conditions (a)--(c), we can find a $\gamma\in \tilde T^0$ such 
that the following intersections are all empty: $d(c^1)\cap 
(d(c^\gamma)\cup\{\beta^\gamma\})$, $A^1\cap A^\gamma$, 
and $B^1\cap B^\gamma$.

Set $p^\gamma_0 = \langle 
c^\gamma\cup\{\langle(\beta^\gamma,\alpha^{\delta_0}),k\rangle\},
\langle
\langle s^\gamma_\alpha,F^\gamma_\alpha\rangle:\alpha\in 
A^\gamma\rangle,\langle\langle t^\gamma_\beta,G^\gamma_\beta 
\rangle:\beta\in B^\gamma\rangle\rangle$, then $p^1\cup 
p^\gamma_0$ is a condition extending both $p^1$ and 
$p^\gamma_0$ that  forces $\dot 
f_{\beta^\gamma}(\alpha^{\delta_0}) = k$. 
This establishes $(*)$.

Let $G^1$ be $\bbb P^1$-generic over $V^0$, and set $V^1 = 
V^0[G^1]=V_2$. 
Let $k_1 = \dot f[G^1](\alpha^{\delta_0})$. By $(*)$, there is a $k\in 
\omega$ and a $\gamma\in \tilde T^0$ so that $k>k_0,k_1$ and 
$\dot f_{\beta^\gamma}[G^1](\alpha^{\delta_0}) = k$. This 
contradicts the 
fact that each $\dot f_{\beta^\gamma}[G^1](\zeta)$ is forced to be 
less than either $k_0$ or $\dot f[G^1](\zeta)$. This establishes the 
claim and the theorem in most cases.

In case $\lambda$ is regular and $\lambda^\omega \le \theta$
fails, we use conditions (2) and (3) of the function $H$ to
carry out a similar construction, replacing (ii) and 
(iv) by
\begin{enumerate}
\item[(ii)$'$] 
$\forall \delta <\gamma\; (\alpha^\gamma > \max\{\alpha^\delta,
\sup ( A(\delta))\})$; and
\item[(iv)$'$] $\sup(H(\beta^\gamma)) = \alpha^\gamma$
(and thus $cf (\alpha^\gamma) = \omega$).
\end{enumerate}
The rest of the argument is very similar to the first
case, and we leave it to the reader to figure out the
details.

In case $cf (\lambda) = \omega$, we write $\lambda = \bigcup_n
\lambda_n$ where $\lambda_n < \lambda_{n+1} < \lambda$ and
the $\lambda_n$ are regular. We   again do a similar
construction, this time producing ordinals $\alpha^\gamma_n
< \lambda_n$; (i), (ii) and (iv) are generalized to
\begin{enumerate}
\item[(i)$''$] if $cf(\gamma) = \omega$, then $\alpha^\gamma_n = 
\sup\{\alpha_n^\delta:\delta<\gamma\}$ for all $n$;
\item[(ii)$''$] if $\gamma$ is a successor ordinal, then 
$\forall n\;\forall \delta <\gamma\; (\alpha_n^\gamma > 
\alpha_n^\delta
$ and, if $A(\delta) \cap \lambda_n$ is bounded in $\lambda_n$,
then  $\alpha_n^\gamma > \sup( A(\delta)))$; and
\item[(iv)$''$] $H(\beta^\gamma)\supseteq
\{ \alpha_n^\delta : n\in\omega\land\delta < \gamma \}$.
\end{enumerate}
The proof continues as before. Notice that
there must be $n \in \omega$ so that $\bar A \cap \{ 
\alpha^\gamma_n
: \gamma < \omega_1\}$ is bounded in $\{ \alpha^\gamma_n :
\gamma < \omega_1 \}$. We complete the argument
with all $\alpha^\gamma$ replaced by $\alpha^\gamma_n$.
\end{Proof}

Suppose that $\lambda<\theta$ are uncountable cardinals that 
satisfy the hypotheses of Theorem~\ref{main-result}. Via the 
translation results in \cite{LL94}, the $(\theta,\lambda)$-good 
set constructed by Theorem~\ref{main-result} gives  a new 
consistent example of a first countable, $<\!\theta$-cwH space $X$ 
that is not $\le\!\theta$-cwH.  The set of non-isolated points of $X$ 
is the union of two disjoint closed discrete sets $D$ and $E$, where 
$|D|= \theta$, $|E| = \lambda$, both $D$ and $E$ are separated, but 
$D$ and $E$ are not contained in disjoint open sets. Thus, $X$ can be 
made cwH by removing the small closed discrete set $E$.

When $\lambda = \theta$ satisfy the hypotheses of 
Theorem~\ref{main-result} and are singular, the space $X$ obtained 
resembles the first countable $<\theta$-cwH not $\le\!\theta$-cwH 
space constructed in~\cite{Fleissner-Shelah}, though the models in 
which the constructions take place may be quite different. In both 
spaces, the set of non-isolated points is the union of two disjoint 
closed discrete sets  of cardinality $\theta$, each of which is 
separated, but that are not contained in a pair of disjoint open sets.

When $\theta$ is singular and greater than $\lambda^\omega$, we 
obtain a first countable space in which cwH fails for the first time 
at $\theta$, yet the space can be made cwH by removing a closed 
discrete set of cardinality $\lambda$.

\section{Characterizing $(\theta,\theta)$-good sets}
We now provide a consistent characterization of the existence of 
$(\theta,\theta)$-good sets in terms of families of integer-valued 
functions. The set-theoretic conditions we require in order to obtain 
this characterization are true in the models obtained by the Levy or 
Mitchell collapse of a large cardinal to $\omega_2$ and when 
$PFA^+$ holds, so this 
characterization may be useful in showing the consistency of ``there 
are no $(\omega_2,\omega_2)$-good sets''.

Recall $\cov(\omega,\theta) = \min\{|\Cal C|: \Cal 
C\subseteq[\theta]^\omega\land \forall B\in 
[\theta]^\omega\;\exists C\in \Cal C\;(B\subseteq C)\}$. It is 
well-known that for all $ n\in\omega$, $\cov(\omega,\omega_n) = 
\omega_n$ and that $\cov(\omega,\theta)>\theta$ for $\theta$ a 
cardinal of uncountable cofinality implies there is an inner 
model with large cardinals.

\begin{Theorem} Assume $\frak b \le \theta$ and 
$\cov(\omega,\theta) \le \theta$.  Then the following are 
equivalent:
\begin{enumerate}
\item[{\em(i)}] there is a $(\theta,\theta)$-good set;
\item[{\em(ii)}] there is an unbounded family $\Cal 
F=\{f_\beta:\beta<\theta\}\subseteq 
\omega^\theta$ such that every $\Cal G\in [\Cal F]^{<\theta}$ is 
bounded and for every $B\in [\theta]^{<\theta}$, $\Cal F\re B = 
\{f_\beta\re B:\beta<\theta\}$ is bounded.
\end{enumerate}
\end{Theorem}
\begin{Proof} (ii) $\Rightarrow$ (i)  This follows from 
Lemma~\ref{functions-imply-good}.

(i) $\Rightarrow$ (ii) Let $A$ be a $(\theta,\theta)$-good set. We 
can assume that each $H_{\beta\alpha}$ (defined in 
Section~\ref{omega-good}) is 
symmetric, i.e., $H_{\beta\alpha}= \{(n,m):(m,n)\in 
H_{\alpha\beta}\}$. Let $\{B_\delta:\delta<\cov(\omega,\theta)\}$ 
enumerate a covering family; for $\delta<\cov(\omega,\theta)$, let 
$r_\delta:\omega\to B_\delta$ be a bijection. Let 
$\{g_\gamma:\gamma\le \frak b\}\subseteq\omega^\omega$ be an 
unbounded family of strictly increasing functions.

For $\beta<\theta$, $\delta<\cov(\omega,\theta)$, and 
$\gamma<\beta$, define a function 
$f_{\beta\delta\gamma}:\theta\to \omega$ by
$$
f_{\beta\delta\gamma}(\alpha) = \left\{ 
\begin{array}{ll}
\max\{n: \exists m\; (g_\gamma(m)\ge r^{-1}_\delta(\alpha) \land 
(\langle\beta,m\rangle,\langle \alpha,n\rangle)\in A)\}& {} \\
0 \quad\mbox{ if the above set is empty or $\alpha\notin 
B_\delta$}&{}
\end{array} \right.
$$
Set $\Cal F = \{f_{\beta\delta\gamma}:\beta<\theta$, 
$\delta<\cov(\omega,\theta)$, and $\gamma<\frak b\}$. We claim
\begin{enumerate}
\item[(I)] $\Cal F$ is unbounded;
\item[(II)] whenever $B\in[\theta]^{<\theta}$, then $\Cal F\re B$ is 
bounded; and
\item[(III)] whenever $\Cal G\in [\Cal F]^{<\theta}$, then $\Cal G$ is 
bounded.
\end{enumerate}

The proofs of (II) and (III) are similar to those for 
Theorem~\ref{char-omega-good}.  For 
example, to prove  (III), fix a $B\in[\theta]^{<\theta}$, and set $A' = 
A\cap((\theta\times \omega)\times(B\times \omega))$.  Then 
$(\infty,\infty)\notin \cl{A'}$, so there is a $g\in \omega^\theta$ so 
that $A'\cap(V_g\times U_g) = \emptyset$.
We claim that $g$ is a bound for $\Cal F\re B$.

To see this, fix $\beta<\theta$, $\delta<\cov(\omega,\theta)$, and 
$\gamma<\frak b$. Take any $\alpha\in B_\delta\less 
r(\{0,1,\dots,g_\gamma(g(\beta))\})$, then $m = r^{-1}(\alpha) 
>g_\gamma(g(\beta))$. Because $g_\gamma$ is increasing, 
$m>g(\beta)$.  Thus, if 
$(\langle\beta,m\rangle,\langle\alpha,n\rangle)\in A'$, we must 
have $n\le g(\alpha)$. Taking the maximum over all such $n$ yields 
$f_{\beta\delta\gamma}(\alpha)\le g(\alpha)$, so clearly 
$f_{\beta\delta\gamma}\le^*g$.  The proof of (II) is similar.

To show that (I) is true, we need the following:

\noindent{\bf Claim:} For every $f\in \omega^\theta$, there is a 
$\beta<\theta$ and a sequence $\{\alpha_m:m\in \omega\}$ so that 
whenever $m\in \omega$, either $f(\alpha_m)\ge m$ and 
$(\langle\beta,m\rangle,\langle\alpha_m,f(\alpha_m)\rangle\in A$, 
or $f(\alpha_m)<m$ and 
$(\langle\beta,m\rangle,\langle\alpha_m,m\rangle)\in A$.

\smallskip
Suppose otherwise. Then there is an $f\in \omega^\theta$ so that for 
all $\beta<\theta$, there is an $m_\beta\in \omega$ such that for 
all $\alpha<\theta$,
\begin{enumerate}
\item[$(*)$] 
$(\langle\beta,m_\beta\rangle,\langle\alpha,m_\beta\rangle)\in A 
\Rightarrow f(\alpha)\ge m_\beta$ and
\item[$(**)$] 
$(\langle\beta,m_\beta\rangle,\langle\alpha,f(\alpha)\rangle)\in 
A\Rightarrow f(\alpha) <m_\beta$.
\end{enumerate}
Define $g\in \omega^\theta$ by $g(\beta) = 
\max\{f(\beta),m_\beta\}+1$. Because $(\infty,\infty)\in \cl{A}$, 
there are $\alpha<\beta<\theta$ with 
$(\langle\beta,g(\beta)\rangle,\langle\alpha,g(\alpha)\rangle)\in 
A$. Now, $H_{\beta\alpha}$ is cdw, so 
$(\langle\beta,m_\beta\rangle,\langle a,f(\alpha)\rangle)\in A$; by 
$(*)$, $f(\alpha)<m_\beta$. By symmetry, $f(\beta)<m_\alpha$. 
Suppose that $m_\beta\le m_\alpha$ (the other case is dual). Then 
$(\langle\beta,m_\beta\rangle,\langle\alpha,m_\beta\rangle)\in 
A$. By $(*)$, $f(\alpha)\ge m_\beta$, a contradiction. This proves 
the claim

\smallskip
To finish the proof of (I), fix $f\in \omega^\theta$, and take 
$\alpha<\theta$ and $\{\alpha_m:m\in \omega\}$ as in the claim.
Find a $\delta<\cov(\omega,\theta)$ so that $\{\alpha_m:m\in 
\omega\}\subseteq B_\delta$.  Define $h\in \omega^\omega$ by 
$h(m) = r^{-1}(\alpha_m)$; then $h$ is finite-to-one. Find a 
$\gamma<\frak b$ so that for infinitely many $m$, $g(m) \le 
g_\beta(m)$. Then for such an $m$, either $f(\alpha) \ge m$ and 
$f_{\beta\delta\gamma}(\alpha_m)<m$, or 
$f_{\beta\delta\gamma}(\alpha_m)\ge m$, so that 
$f_{\beta\delta\gamma}\not\le^+f$.\end{Proof}

\section{Questions} We have shown for many cardinals $\lambda\le 
\theta$ that ``there is a $(\theta,\lambda)$-good set'' is consistent.
On the other hand, under $GCH$, there are no 
$(\theta,\lambda)$-good sets when $\omega_1\le cf(\lambda)\le 
\lambda<\theta$, so for these cardinals, the existence of a 
$(\theta,\lambda)$-good set is independent of $ZFC$. There are no 
$(\theta,\theta)$-good sets when $\theta$ 
is singular of countable cofinality (see \cite{LL94}), so we ask:

\begin{enumerate}
\item[(1)] Is it consistent to have a 
$(\theta,\lambda)$-good set when $cf(\theta) = \omega$ and 
$\omega_1\le cf(\lambda)\le\lambda<\theta$?
\item[(2)] Suppose that $\omega = cf(\lambda)<\lambda$. Is there, in 
$ZFC$, a cardinal $\theta$ such that $\lambda<\theta\le 
\lambda^\omega$ and there is a $(\theta,\lambda)$-good set?
\end{enumerate}

Of course, the most important question, originally asked by  Dow and
Todor\v{c}evi\'c, is:
\begin{enumerate}
\item[(3)]  Does $ZFC$ imply the existence of an 
$(\omega_2,\omega_2)$-good set?
\end{enumerate}
Todor\v{c}evi\'c \cite{Todorcevic-fan} showed that 
$\square(\omega_2)$ 
implies that there is a $(\omega_2,\omega_2)$-good set. So at least 
a weakly compact cardinal is required to produce a model with no 
$(\omega_2,\omega_2)$-good sets. Fleissner used 
$E^\omega_{\omega_2}$ (i.e., ``there is a non-reflecting stationary 
subset of $\omega_2$ consisting of ordinals of countable 
cofinality'') to construct a first countable, $<\!\omega_2$-cwH 
space that is not $\le\!\omega_2$-cwH, so $E^\omega_{\omega_2}$ 
can also be used to produce an $(\omega_2,\omega_2)$-good set.

Recall that Beaudoin (and independently, Magidor) showed that $PFA$ 
is consistent with $E^\omega_{\omega_2}$, while 
$PFA^+$ implies that stationary sets reflect. We conclude with:
\begin{enumerate}
\item[(4)] Does $PFA^+$ imply that there are no 
$(\omega_2,\omega_2)$-good sets?
\end{enumerate}

\begin{tabbing}
\noindent{\bf Addresses:} \=
{\sc Department of Mathematics, University of Tuebingen, 
Tuebingen,}\\
{}\>{\sc Germany}\\
\\
{}\>{\sc Department of Mathematics, Union College, Schenectady, 
NY,}\\
{}\>{\sc  12308, USA}
\end{tabbing}

\begin{tabbing}
\noindent{\bf E-mail: }
\= \verb+jobr@michelangelo.mathematik.uni-tuebingen.de+\\
{}\>\verb+laberget@gar.union.edu+
\end{tabbing}


\begin{thebibliography}{FS}

\bibitem[Br]{Br}
J. Brendle.
\newblock Notes.

\bibitem[FS]{Fleissner-Shelah}
W.~G. Fleissner and S. Shelah.
\newblock Incompactness at singulars.
\newblock {\em Topology Appl.}, 31:101--107, 1989.

\bibitem[G]{G}
G. Gruenhage.
\newblock $k$-spaces and products of closed images of metric 
spaces.
\newblock {\em Proc. Amer. Math. Soc.}, 80(3):477--482, 1980.

\bibitem[H]{Hausdorff}
F. Hausdorff.
\newblock Die Graduierung nach dem Endverlauf.
\newblock Abhandlungen der K\"oniglich S\"achsischen Gesellschaft 
der Wissenschaften.
\newblock {\em Mathematisch--Physische Klasse\/} 31(1909), 
296--334.

\bibitem[J]{Judah}
H. Judah.
\newblock Private communication.

\bibitem[K]{Koszmidor-Kurepa}
P. Koszmider.
\newblock Kurepa trees and topological non-reflection.
\newblock Preprint.

\bibitem[Ku]{Kunen-diss}
K. Kunen.
\newblock Inaccessibility properties of cardinals.
\newblock Ph.D. dissertation, Stanford, 1968.



\bibitem[LL]{LL94}
T. LaBerge and A. Landver.
\newblock Tightness in products of fans and psuedo-fans.
\newblock Preprint.

\bibitem[R]{Rothberger}
F. Rothberger.
\newblock Sur les familles ind\'enombrables de suites de nombres 
naturels et les probl\`emes concernant la propri\'et\'e $C$.
\newblock {\em Proc. Camb. Phil. Soc.}, 37:109--126, 1941.

\bibitem[T]{Todorcevic-fan}
S. Todor\v{c}evi\'c.
\newblock My new fan.
\newblock Preprint.

\end{thebibliography}
\end{document}